# INVARIANT BAYESIAN ESTIMATION ON MANIFOLDS

### By Ian H. Jermyn

### *INRIA*

A frequent and well-founded criticism of the maximum a posteriori (MAP) and minimum mean squared error (MMSE) estimates of a continuous parameter $\gamma$ taking values in a differentiable manifold $\Gamma$ is that they are not invariant to arbitrary "reparameterizations" of $\Gamma$. This paper clarifies the issues surrounding this problem, by pointing out the difference between coordinate invariance, which is a sine qua non for a mathematically well-defined problem, and diffeomorphism invariance, which is a substantial issue, and then provides a solution. We first show that the presence of a metric structure on $\Gamma$ can be used to define coordinate-invariant MAP and MMSE estimates, and we argue that this is the natural way to proceed. We then discuss the choice of a metric structure on $\Gamma$. By imposing an invariance criterion natural within a Bayesian framework, we show that this choice is essentially unique. It does not necessarily correspond to a choice of coordinates. In cases of complete prior ignorance, when Jeffreys' prior is used, the invariant MAP estimate reduces to the maximum likelihood estimate. The invariant MAP estimate coincides with the minimum message length (MML) estimate, but no discretization or approximation is used in its derivation.

**1. Introduction.** Statistical estimation is a very old field, but despite that many questions remain unanswered and debates about the best way to proceed are plentiful. From a probabilistic point of view, all the information about a quantity of interest taking values in a space $\Gamma$ is contained in a probability measure on $\Gamma$. If it is deemed necessary to single out a particular point $\gamma \in \Gamma$ for some purpose, a loss function $L : \Gamma \times \Gamma \to \mathbb{R} : (\gamma, \gamma') \rightsquigarrow L(\gamma, \gamma')$ is defined describing the cost inherent in taking the true value of the quantity to be $\gamma$ when it is in fact $\gamma'$. The mean value of the loss as a function of $\gamma$ can be computed using the probability measure, whereupon one can, for example, choose that point $\hat{\gamma} \in \Gamma$ that minimizes the mean loss as one's estimate of the true value of $\gamma$.

In some cases, especially those closely linked to a specific application, the loss function will be fully dictated by circumstance. In this case, the invariance issues discussed in this paper do not arise. However, in many other cases, and for the purposes of theoretical analysis, estimates are needed in the absence of any clear knowledge of what the real loss is. Indeed, there may not be a "real loss." In these cases, generic loss functions are required, and indeed are currently widely used, in both theory and practice. These generic loss functions should satisfy two criteria: they must be well-defined, and they must not introduce implicit bias that is not present in the models. The latter is best expressed by saying that in the absence of prior knowledge about the loss function, the loss function should not introduce prior knowledge about the parameters to be estimated. This is an application of the principle that if two people have the same knowledge, then they should make the same inferences.

In the case that $\Gamma$ is a differentiable manifold, difficulties arise. Two popular choices of generic loss function are the negative of a delta function and the squared difference of coordinates, leading to maximum a posteriori (MAP) and minimum mean squared error (MMSE) estimates, respectively. In order for these quantities to be well defined, two things are necessary: an underlying measure in order to define the delta function loss, and a distance function in order to define the squared error. The existence of these quantities is normally ignored, or equivalently they are assumed to take on particular forms. The resulting loss functions are not coordinate-invariant, and hence are ill-defined in general coordinate systems, thus violating the first criterion. This lack of coordinate invariance leads to the paradox that two people with the same knowledge can construct different estimates simply by choosing to use different coordinate systems, for example, polar rather than rectangular. Even if the definitions are made coordinate-invariant, and hence well-defined, the resulting loss functions still violate the second criterion in general. The estimates are not invariant to diffeomorphisms, which "mix up" the points of $\Gamma$ ("reparametrizations"), and therefore necessarily introduce extra information about these points.

The purpose of this paper is to correct the above situation. We define compatible, coordinate-invariant MAP and MMSE estimates by introducing a Riemannian metric on $\Gamma$, and argue that this is the natural way to achieve such invariance. This satisfies the first criterion. The introduction of a metric raises the question of how to choose this extra structure, and we argue that in the case of Bayesian estimation, imposing the second criterion renders the choice of metric unique.

The main results of the paper as regards Bayesian estimation are the following:

(a) The metric on $\Gamma$ should be the pullback by the model function of the natural metric on every measure space.

(b) Invariant MAP estimates should be defined using the density of the posterior measure with respect to the measure derived from this metric.

(c) Invariant MMSE estimates should be defined by using, in place of the squared error, the squared geodesic distance based on this metric.

(d) In conditions of "complete ignorance," that is, conditions in which the prior probability measure is Jeffreys' prior, MAP estimates always reduce to maximum likelihood (ML) estimates, in contrast to much Bayesian argument and practice.

(e) The invariant MAP estimate coincides with the MML estimate described by Wallace and Freeman (1987), except that no discretization of $\Gamma$ is required and no approximations are made.

The rest of the paper is structured thus. In Section 2 we discuss the failure of invariance for MAP estimation on manifolds and its causes. In Section 3 we describe how both this problem, and the related failure of invariance for MMSE estimates, can be solved by endowing the manifold with a metric structure, and we argue that this is the natural solution to the problem. In Section 4 we discuss the choice of metric structure, and use a simple invariance argument to render this choice unique. In Section 5 we discuss the conclusions of the report and related work.

The material on the differential geometry of measure spaces and its connection to Jeffreys' prior may be known to geometrically minded statisticians. We include it here for completeness, and to emphasize its coordinate-invariant nature.

**2. The problem.** To illustrate the problem, we examine the maximization of a probability density function (p.d.f.) on a manifold of dimension $m$. Let the manifold be $\Gamma$, a point in $\Gamma$ being denoted $\gamma$. We are given a probability measure $\mathbf{Q}$ on $\Gamma$, which we may view as the posterior in an MAP estimation task, although this is not important at this stage. We are also given two systems of coordinates on $\Gamma$, $\theta\colon\Gamma\to\mathbb{R}^m$ and $\phi\colon\Gamma\to\mathbb{R}^m$. (We ignore questions of topology that might require us to use more than one coordinate patch; the issue is not central to the discussion here.)

Expressed in terms of the first set of coordinates $\theta$, and the corresponding measure $d^m\theta(\gamma)$ on $\Gamma$, we find $\mathbf{Q} = Q_\theta(\theta(\gamma))\,d^m\theta(\gamma)$, where $Q_\theta(\theta(\gamma))$ is a function. We now separate the function $Q_\theta$ from the measure and find the argument of its maximum value $\theta_{\max}\in\mathbb{R}^m$, giving an estimate of $\gamma$, $\hat{\gamma}_\theta = \theta^{-1}(\theta_{\max})$.

We may choose to express $\mathbf{Q}$ in another coordinate system, $\phi\colon\Gamma\to\mathbb{R}^m$. Using the measure defined by this coordinate system, we find that $\mathbf{Q} = Q_\phi(\phi(\gamma))\times d^m\phi(\gamma)$. If we now follow the same procedure as before, and find the argument of the maximum value of $Q_\phi$, $\phi_{\max}$, we find another estimate, $\hat{\gamma}_\phi = \phi^{-1}(\phi_{\max})$.

The problem is the following. Suppose that the two coordinate systems are related by a function $\alpha \colon \mathbb{R}^m \to \mathbb{R}^m$, so that $\theta(\gamma) = \alpha(\phi(\gamma))$. In this case, the measures with respect to the two coordinate systems are related by $d^m\theta(\gamma) = J[\alpha](\phi(\gamma))\, d^m\phi(\gamma)$, where $J[\alpha](\phi(\gamma))$ is the Jacobian of the coordinate transformation. This in turn means that the functions $Q_\theta$ and $Q_\phi$ are related by $Q_\phi(\phi(\gamma)) = Q_\theta(\theta(\gamma))J[\alpha](\alpha^{-1}(\theta(\gamma)))$.

The consequence is that the estimates obtained by maximizing $Q_\theta$ and $Q_\phi$ are different, due to the presence of the Jacobian factor. Apparently our estimate of $\gamma$ depends upon the choice of coordinates, or in effect upon the whim of the person making the estimate. This may seem surprising: one thinks of the question "What is the most probable point in $\Gamma$?" and, by analogy with the discrete case, one expects an invariant answer.

The difference between the continuous and the discrete cases means, however, that the question being asked in the continuous case is not the previously cited one at all, but a slightly more complicated version. Given a coordinate system, $\theta$, the question being asked is, "What is the infinitesimal volume element $\theta^{-1}(dz)$ in $\Gamma$ (where $dz$ is an infinitesimal coordinate volume in $\mathbb{R}^m$) that is most likely to contain the true point in $\Gamma$?" (We use the notation $f^{-1}$ both for the inverse of a map $f \colon A \to B$, $f^{-1} \colon B \to A$, and for the pullback $f^{-1} \colon 2^B \to 2^A \colon B \supset Y \rightsquigarrow \{a \in A \colon f(a) \in Y\}$. Context serves to distinguish the two usages.) Using a different coordinate system, $\phi$ on the other hand, the question is "What is the infinitesimal volume element $\phi^{-1}(dz)$ that is most likely to contain the true point in $\Gamma$?" In general, $\theta^{-1}(dz) \neq \phi^{-1}(dz)$. It is then clear that different answers are to be expected using different coordinate systems, because the question being asked is different in each case.

A simple example of the above is provided by a Gaussian measure in two dimensions with zero mean and covariance the identity. This measure can be expressed in rectangular or in polar coordinates:

$$\Pr(\vec{x}) = dx\, dy\, Z^{-1} e^{-(x^2+y^2)} = dr\, d\theta\, Z^{-1} r e^{-r^2}.$$

In the first case, the maximum density procedure leads to $\hat{x} = \hat{y} = 0$, while in the second it leads to $\hat{r} = 1/\sqrt{2}$ and an indeterminate value for $\hat{\theta}$. In this simple case, one can see the error clearly, but in more complex or less intuitive cases the same phenomenon arises and passes unnoticed.

From a measure-theoretic point of view, what is happening is clear. The functions $Q_\theta$ and $Q_\phi$ are probability density functions. Any p.d.f. is defined with respect to an underlying measure. The Radon–Nikodym derivative of the probability measure with respect to the underlying measure then gives the p.d.f. In the scenario just described, two different underlying measures are being used: $d^m\theta(\gamma)$ and $d^m\phi(\gamma)$. To expect them to yield the same results is unreasonable.

If one concentrates on the underlying measure, then there is no problem. In terms of $\theta$, the underlying measure is $d^m\theta(\gamma)$, while in terms of $\phi$, the same underlying measure is $J[\alpha](\phi(\gamma))\,d^m\phi(\gamma)$. Integration of either of these over a fixed subset of $\Gamma$ will produce the same result: they are the same measure. Using this fixed measure, the problem disappears: in terms of $\phi$, the p.d.f. with respect to the underlying measure is $Q_\theta(\alpha(\phi(\gamma))) = Q_\theta(\theta(\gamma))$. The maxima of $Q_\theta(\alpha(\phi(\gamma)))$ with respect to $\phi$ agree completely with those of $Q_\theta(\theta(\gamma))$ with respect to $\theta$, in the sense that $\theta_{\max} = \alpha(\phi_{\max})$, which implies that $\theta^{-1}(\theta_{\max}) = \phi^{-1}(\phi_{\max})$. The points in $\Gamma$ that we find are the same. The problem is that, given an arbitrary coordinate system, we do not know which choice of coordinate is "correct," and hence what the estimate should be. By effectively focusing on measures on $\mathbb{R}^m$, the coordinate space, rather than on underlying measures on $\Gamma$, the problem is created. How then to define, in a coordinate-invariant way, an underlying measure with respect to which to take the Radon–Nikodym derivative?

A similar situation arises with respect to MMSE estimates, which also lack invariance under general changes of coordinates. It is equally true that the mean itself has no coordinate-invariant meaning, and for the same reasons. In calculating both the error and the mean, one is faced with adding or subtracting certain values. If these operations are performed on the coordinate values in a particular coordinate system, they will change with a change of coordinates. Equally, one cannot add or subtract points of $\Gamma$ directly; such operations are not defined unless $\Gamma$ possesses an algebraic structure of some kind, for example, is a vector space.

In practice, what is crucial to the MMSE estimate is the notion of a distance between two points in $\Gamma$. If a global Euclidean coordinate system exists, this is given by the squared error, but in general this is not the case. If we wish to consider MMSE estimates in general coordinate systems, we must be able to define distances in a coordinate-invariant manner.

**3. Coordinate-invariant estimates.** If one wishes to discuss measures and distances using an arbitrary set of coordinates, one must express the mathematics in a way that allows for this eventuality. Not to do so means that symbols such as $d^m\theta$ are not defined. The natural way to express both geometric and measure-theoretic information about manifolds in a way that is manifestly free of coordinates, but that nevertheless allows the derivation of an expression in terms of an arbitrary coordinate system with the greatest of ease, is the language of forms. Readers not familiar with this language may wish to look at the Appendix, where we provide a brief introduction to forms and their uses, or at the book by Choquet-Bruhat, DeWitt-Morette and Dillard-Bleick (1977).

We are interested in probability measures. These can be integrated over $m$-chains, for example, the whole manifold $\Gamma$, and as such are $m$-forms. In addition, they must be positive and normalized, so that they are probability $m$-forms. The answer to the first of the questions at the end of the last section is then: define an $m$-form, since these are, by definition, coordinate-invariant. The answer to the second question would seem to be: define a distance function. In practice, the following considerations push us strongly in one direction: the introduction of a Riemannian metric on the manifold $\Gamma$.

First, the introduction of a metric allows us simultaneously to answer both of the questions posed at the end of the last section. Starting from the metric, we can derive an $m$-form and use this as the underlying measure. We can also define a distance function, as the geodesic distance between two points.

Second, if we are to introduce notions both of "volume" (via an underlying measure) and of "length" (via a distance function), it is sensible that these notions be compatible. Otherwise there is no reason to believe that the resulting estimates will bear any relation to one another. The use of a metric to define both the underlying measure and the distance function ensures that maps that preserve lengths preserve volumes also, or, even more intuitively, that the volume of a small cube is given by the product of the lengths of its sides.

The final consideration is intuition in practice. Manifolds with a measure but no metric are strange objects. They do not correspond to our intuition of a surface or volume at all. The space of volume-preserving diffeomorphisms is much larger than the space of isometries, and allows severe distortions. An example is the mixing of two incompressible immiscible fluids. The initial "drop of oil in water" may end up smoothly distorted into dramatically different shapes. The parameter spaces that we consider intuitively possess "metric-like" properties, even if these are not well defined. For a one-dimensional $\Gamma$, for example, the numbers that represent different parameter values indicate something more than the topological, although a precise interpretation may not be available. If we wish to be able to describe these geometric properties of the manifold as well as its measure-theoretic properties, a metric is necessary. In addition, it is quite hard to write down an expression for a measure on a manifold without implicitly assuming a metric. In practice, this means that metrics appear, albeit disguised, in the expressions for many probability measures. Gaussian measures are one example, where an inner product is used to define the exponent. An inner product on a vector space is equivalent to a constant metric, which allows identification of each tangent space with the vector space itself. In many other cases, the assumption of a Euclidean metric is made manifest by the appearance of an orthogonal inner product.

What then is a Riemannian metric and how does it define a measure? A metric $\mathbf{h}$ is the assignment, to each point $\gamma$ of $\Gamma$, of an inner product on the tangent space $T_\gamma\Gamma$ at $\gamma$. This is detailed in the Appendix, where it is further explained how the existence of a metric allows us to map functions to $m$-forms using the Hodge star. Given a function $f$ on $\Gamma$, we can thus create an $m$-form, that is, a measure $\star_\mathbf{h} f$. The choice of function $f$ is dictated by compatibility between the measure-theoretic and geometric aspects of the manifold. By choosing $f$ to be $\mathbb{I}$, the function identically equal to 1, the resulting $m$-form is preserved by isometries; in other words, maps that preserve length preserve volume also.

Being a form, the quantity $\mathbf{U_h} = \star_\mathbf{h}\mathbb{I}$ is invariantly defined. This is clear first because no coordinate system was used in its construction, but it can also be verified in detail. As described in the Appendix, the expression for this form in the coordinate basis of coordinates $\theta$ is

$$\mathbf{U_h} = \star_\mathbf{h}\mathbb{I} = |\mathbf{h}|_\theta^{1/2}\, d^m\theta,$$

where $|\mathbf{h}|_\theta$ is the determinant of the metric components in the $\theta$ coordinate basis, and $d^m\theta$ is the coordinate basis element for the space of $m$-forms. To see the invariance of this measure explicitly, note that a change of coordinates $\alpha$ introduces a factor of $J[\alpha](\phi(\gamma))$ from $d^m\theta$, while the transformation of the determinant of the metric matrix elements from one basis to another introduces a factor of $J[\alpha](\phi(\gamma))^{-1}$. Thus, expressed in any coordinate system, the form of the measure is identical: $|\mathbf{h}|_\theta^{1/2}\, d^m\theta = |\mathbf{h}|_\phi^{1/2}\, d^m\phi$. To stress the point once again: the measure $d^m\theta(\gamma)$ has no coordinate-invariant meaning. If we try to express a measure in a general coordinate system in this way, we literally do not know what we are talking about.

3.1. *Maximum density estimates.* Given a probability $m$-form $\mathbf{Q}$, and another positive $m$-form $\mathbf{U}$, one defines the p.d.f. of $\mathbf{Q}$ with respect to $\mathbf{U}$ by division:

$$(3.1) \qquad Q = \frac{\mathbf{Q}}{\mathbf{U}}.$$

This is the equivalent of the Radon–Nikodym derivative in the language of forms. What now becomes of maximum density estimation? We simply have to use $\mathbf{U_h}$ in (3.1). If we choose a particular coordinate system $\theta$, so that $\mathbf{Q} = Q_\theta\, d^m\theta$ and $\mathbf{U_h} = |\mathbf{h}|_\theta^{1/2}\, d^m\theta$, then we have

$$(3.2) \qquad Q = |\mathbf{h}|_\theta^{-1/2} Q_\theta.$$

The left-hand side of this equation is invariant to changes in coordinates. These will produce equal Jacobian factors in both the numerator and the denominator of (3.2), which will thus cancel out. Note also that this p.d.f.

does not result simply from a choice of coordinates. Although it may be possible to find a system of coordinates in which the determinant of the metric is constant, this is misleading in two ways. First, what is really happening is that a metric is being chosen. The naive approach really means choosing a metric whose determinant is constant in the coordinate system you already have, which is not a coordinate-invariant procedure. Second, in more than one dimension, although the determinant of the metric may be constant, it may not be possible to find a system of coordinates in which the metric itself is constant. This would imply that the manifold is flat, a statement that is coordinate-invariant and may not be true.

3.1.1. *Expression in terms of a delta function loss.* Usually the maximum density estimate is regarded as derived from the use of a particular loss function, $\delta(\theta(\gamma), \theta(\gamma'))$ on $\Gamma$. Given a probability $m$-form expressed in terms of $\theta$, $Q_\theta(\theta) \, d^m\theta$, this leads to the familiar recipe $\hat{\gamma}_\theta = \theta^{-1}(\arg\max_\theta Q_\theta(\theta))$, in apparent contradiction to the previous discussion. From this point of view, there is no need to define a p.d.f. at all, since we were merely integrating with respect to the probability measure. What is going on?

The answer of course involves the same concepts as above. The quantity $\delta(\theta(\gamma), \theta(\gamma'))$ is not invariantly defined, since the measure against which to integrate it has not been given. In our context, the delta function (in fact there are effectively $m$ of them) is best viewed as the identity map from $\Lambda^p\Gamma$, the space of $p$-forms on $\Gamma$, to itself. As such, it is a $p$-form at its first argument (a point in $\Gamma$) and an $(m-p)$-form at its second argument (another point in $\Gamma$). It can thus be integrated against a $p$-form to produce another $p$-form. When $p = 0$, we recover the usual delta function that evaluates a function at its first argument. In our case, however, we wish to integrate the delta function against an $m$-form, and thus $p = m$. The delta function is thus an $m$-form at its first argument and a 0-form, or function, at its second argument. The result of integrating it against the posterior measure is thus an $m$-form, and to create a function that we can maximize, we need to use the Hodge star. This again introduces the factor of $|\mathbf{h}|_\theta^{-1/2}$ that we see in (3.2) and that is implicit in (3.1).

An alternative point of view is to consider the delta function as a map from $\Lambda^p\Gamma$ to $\Lambda^{(m-p)}\Gamma$, making it an $(m-p)$-form at its first argument and a $p$-form at its second argument. In order to integrate this against a $p$-form, we can use the inner product on $\Lambda^p$ described in (A.2) of the Appendix. In our case, this point of view makes the delta function a 0-form (function) at its first argument and an $m$-form at its second. The result of the integration is thus a function as required for maximization, but now we find that the use of the inner product has already introduced the factor of $|\mathbf{h}|_\theta^{-1/2}$, thus giving the same result as in the other two methods.

There is thus no conflict between these different ways of speaking.

3.2. *MMSE estimates.* Suppose we are given a distance function. That is, we are given a symmetric map $d : \Gamma \times \Gamma \to \mathbb{R}^+$, obeying the triangle inequality and such that $d(\gamma, \gamma) = 0$. Given a point $\gamma$, we define the function

$$d_\gamma(\gamma') = d(\gamma, \gamma').$$

We can now define the coordinate-invariant form of the mean squared error, which we will call the *mean squared distance*, as

$$(3.3) \qquad L(\gamma) = \int_\Gamma (d_\gamma)^2 \mathbf{Q},$$

where $\mathbf{Q}$ is as usual a probability $m$-form. In terms of a particular coordinate system $\theta$ on $\Gamma$, one has

$$L(\theta) = \int_{\theta(\Gamma) \subset \mathbb{R}^m} d^m \theta' \, Q_\theta(\theta') \, d_\theta^2(\theta, \theta'),$$

where $d_\theta$ is the expression for the length in terms of the given coordinates.

Having defined the mean squared distance $L$, we can now define the minimum mean squared distance (MMSD) estimate as the set of minimizers of $L(\gamma)$.

All that remains is to use the metric to define a distance function that we can use in (3.3). Below we recap this material from differential geometry, phrasing it in a manifestly coordinate-invariant way, and emphasizing the difference between coordinate invariance and invariance to diffeomorphisms, which is a coordinate-invariant and therefore content-full concept. We first define the notion of the length of a path, and then define the distance between two points as the length of a minimum length path between them.

Let $I$ be an interval of the real line, considered as a manifold (i.e., without the structure of a field). Let $p_0$ and $p_1$ be the elements of its boundary. Let $\pi : I \to \Gamma$ be an embedding of $I$ in $\Gamma$ such that $\pi(p_0) = \gamma$ and $\pi(p_1) = \gamma'$. To define the length of the path (i.e., its volume), we need a 1-form on $I$, or in other words a measure, which we will then integrate over $I$. Now, however, we have an invariance criterion: we must ensure that the length we calculate depends only on the image of $I$ in $\Gamma$, and not on the precise mapping of points of $I$ to points of $\Gamma$. This amounts to saying that replacing $\pi$ by $\pi\varepsilon$, where $\varepsilon$ is an arbitrary boundary-preserving diffeomorphism, should not change the resulting length. Note that unlike coordinate invariance on $I$, which follows as soon as we integrate over the coordinates, this condition is a substantive one. As argued in the Appendix, the only way to ensure this is to construct a metric on $I$ by pulling back a metric from $\Gamma$, and then using this metric in the normal way to construct a 1-form. We thus pull back the metric $\mathbf{h}$ on $\Gamma$ to give a metric $\pi^*\mathbf{h}$ on $I$. We then use the Hodge star of this metric to map $\mathbb{I}$ to a 1-form that can be integrated on $I$. In notation,

$$(3.4) \qquad l(\pi) = \int_I \star_{\pi^*\mathbf{h}} \mathbb{I}.$$

To illustrate the ability to derive an expression in an arbitrary coordinate system from the coordinate-invariant expression (3.4), we introduce a coordinate system $t: I \to \mathbb{R}$ on $I$ with a corresponding coordinate basis given by $\frac{\partial}{\partial t}(p)$, and a coordinate system $\theta$ on $\Gamma$ with a corresponding coordinate basis given by $\frac{\partial}{\partial \theta^i}(\gamma)$. In these bases, the (single) component of the pulled back metric can be found to be

$$(\pi^* \mathbf{h})_p \left( \frac{\partial}{\partial t}(p), \frac{\partial}{\partial t}(p) \right) = \mathbf{h}_{\pi(p)} \left( \frac{d\pi^i}{dt}(p) \frac{\partial}{\partial \theta^i}(\pi(p)), \frac{d\pi^j}{dt}(p) \frac{\partial}{\partial \theta^j}(\pi(p)) \right)$$

$$= h_{\pi(p),ij} \frac{d\pi^i}{dt}(p) \frac{d\pi^j}{dt}(p),$$

where $h_{ij}$ are the components of the metric $\mathbf{h}$ in the $\theta$ coordinate system. Thus the result is simply the length of the tangent vector to the path $\pi$ in the metric $\mathbf{h}$. Rewriting (3.4) in terms of this expression, we find that

$$l(\pi) = \int_a^b dt \left( h_{\pi(t),ij} \frac{d\pi^i}{dt}(t) \frac{d\pi^j}{dt}(t) \right)^{1/2},$$

where we have abused notation by using the same symbol $\pi$ for the map from $I$ to $\Gamma$ and its expression in terms of coordinates. The points $a \in \mathbb{R}$ and $b \in \mathbb{R}$ are the coordinate values of $p_0$ and $p_1$, respectively.

Given the length of a path, we can now define the distance between two points as

$$d_\gamma(\gamma') = d(\gamma, \gamma') = \min_{\pi \in \Pi(\gamma, \gamma')} l(\pi),$$

where $\Pi(\gamma, \gamma')$ is the space of paths with endpoints $\gamma$ and $\gamma'$. This distance is coordinate-invariant, and can be used in (3.3). For a general metric it is of course hard to derive an analytic expression for $d$.

In the case that the metric is Euclidean, $L$ reduces to the mean squared error, as it should. The resulting MMSD estimate is then the mean, that is, the MMSE estimate, and is unique. In other cases, the MMSD estimate provides a generalized mean, known as the "Karcher mean," first introduced by Karcher (1977) as the center of mass on a Riemannian manifold. It is a set of points in $\Gamma$, each of which minimizes the mean squared distance to every other point of $\Gamma$. Note that the set of minimizers may contain more than one point of $\Gamma$. This does not present a problem as such. It simply means that from the point of view of the mean squared distance loss function these points are equivalent.

**4. Bayesian estimation and the choice of metric.** We have argued that in order to define coordinate-invariant and consistent maximum density and MMSE estimates, one should use a metric on the manifold $\Gamma$. We now turn to the question that we have been conspicuously avoiding. How is one to choose a metric on $\Gamma$?

Thus far, we have been dealing solely with a manifold $\Gamma$ and a probability measure $\mathbf{Q}$ on this manifold. In this abstract situation, it seems that the above question has no good answer, which is unsurprising. We turn now, however, to the case that is usually of interest: when $\mathbf{Q}$ is a posterior probability measure derived from a model function and a prior using Bayes' theorem.

We introduce the data space, $X$. We assume that this has sufficient structure to allow the following constructions, and in practice it can be supposed to be either a countable set or a manifold. On $X$ one can define the space of measures $\mathcal{M}(X)$. The space of probability measures, $\mathcal{S}(X)$, is a proper subset of the cone of positive measures. This set has a complicated boundary even in the case where $X$ is countable; when $X$ is not countable, there are also measures with singular components, which complicate things still further. We avoid these difficulties by assuming that all measures with which we will deal lie in the interior of $\mathcal{M}(X)$ and, where appropriate, are nonsingular.

We are free to choose coordinates on $\mathcal{M}(X)$ as on any manifold. One choice is to describe measures as $n$-forms, in which case the space $\mathcal{S}(X)$ becomes the space of probability $n$-forms. A model function is a map $\Lambda : \Gamma \to \mathcal{M}(X)$ associating to each point $\gamma \in \Gamma$ a (probability) measure on $X$. We will assume that this map is a regular embedding, so that the image of $\Gamma$ with the differentiable structure induced by $\Lambda$ is a submanifold of $\mathcal{S}(X)$.

### 4.1. *An invariance criterion.*

We now use this extra structure, which is present in any real estimation problem, to argue for a unique choice of metric on $\Gamma$. The argument rests on one simple idea: that all information about the parameters not contained in the data be contained in the prior measure, or in other words, that all information that distinguishes one point of $\Gamma$ from another should come either from their correspondences with probability measures on $X$ (condition 1) or from the prior measure on $\Gamma$ (condition 2). It is the probability measures on $X$ alone that determine the relationship between the points in $\Gamma$ and the observations represented by points in $X$, and the way that these measures are parameterized serves to determine the meaning of the points in $\Gamma$ and not the other way around. Any other information in addition to the data we have at hand should be described by the prior. Any metric that we choose on $\Gamma$ should respect this principle, and not introduce any extra information about points in $\Gamma$. This is the second criterion.

The fact that it is not the identity of individual points in $\Gamma$ that is important, but merely their correspondence with probability measures on $X$, means that it is only the image of $\Gamma$ in $\mathcal{M}(X)$ that counts. This image is invariant under the replacement of $\Lambda$ by $\Lambda\varepsilon$, where $\varepsilon : \Gamma \to \Gamma$ is a diffeomorphism. A model function is thus an equivalence class of maps $\{\Lambda\varepsilon\}$. The conclusion from condition 1 is thus that inference should be invariant under the replacement of $\Lambda$ by $\Lambda\varepsilon$, where invariant means that the image of the estimate by the model function is preserved. This diffeomorphism invariance, although superficially similar to a change of coordinates, is defined independently of any change of coordinates, and as such is a substantive restriction.

There are only two ways to achieve this aim. One is to pick a particular representative of the equivalence class of maps $\{\Lambda\varepsilon\}$ and to define a metric on the corresponding copy of $\Gamma$. This metric can then be pulled back to other members of the equivalence class using the maps $\varepsilon$. Although this will satisfy condition 1, the selection of a particular member of the equivalence class to be endowed with a particular metric implies that we already know something about the points in $\Gamma$ independently of their correspondence with probability measures on $X$. Otherwise, how could we know to which points of $\Gamma$ to assign which values of the metric? This is exactly the type of information that should be included in the prior, and thus the procedure described in this paragraph violates condition 2.

The second approach is to pull back a metric from $\mathcal{M}(X)$ to each equivalent copy of $\Gamma$ using $\Lambda\varepsilon$. [Since an embedding is a full rank immersion, the pulled back metric will be a proper Riemannian structure on $\Gamma$ if $\mathcal{M}(X)$ is a proper Riemannian manifold.] Such a metric automatically satisfies the consistency conditions introduced by the maps $\varepsilon$ between members of the equivalence class: $\Lambda\varepsilon^*\mathbf{g} = \varepsilon^*\Lambda^*\mathbf{g}$, where $\mathbf{g}$ is a metric on $\mathcal{M}(X)$, and thus our results will depend solely on the image of $\Gamma$ in $\mathcal{M}(X)$. In addition, we were not required to pick a particular member of the class a priori, since each member of the equivalence class gets its own consistent metric induced by its own model function. Thus both condition 1 and condition 2 are satisfied.

We are thus in a position to define a metric and underlying $m$-form on $\Gamma$ that satisfies the invariance criterion stated at the beginning of this section by pulling back a metric from $\mathcal{M}(X)$. We lack only one thing: a metric on $\mathcal{M}(X)$ to pull back.

4.2. *Metrics on* $\mathcal{M}(X)$. The first thing we must do is to define what we mean by the tangent space to $\mathcal{M}(X)$. Since we are using $n$-forms as coordinates on $\mathcal{M}(X)$, and since the space of signed measures is linear, it is easy to see that a tangent vector to $\mathcal{M}(X)$ can be identified with an $n$-form. If we restrict attention to $\mathcal{S}(X)$, this $n$-form must integrate to zero to preserve normalization. Then, at a point $\mathbf{T} \in \mathcal{M}(X)$, an inner product between two tangent vectors $\mathbf{v}_1$ and $\mathbf{v}_2$ is given by

$$\mathbf{g}(\mathbf{v}_1, \mathbf{v}_2) = \int_X \mathbf{T} \frac{\mathbf{v}_1}{\mathbf{T}} \frac{\mathbf{v}_2}{\mathbf{T}},$$

where we have identified the abstract tangent vectors $\mathbf{v}$ with their expression as $n$-forms. Note that the divisions are well defined because $\mathbf{T}$ is positive. The justifications for this choice as the only reasonable metric on $\mathcal{M}(X)$ are many, and we do not reiterate them here. Interested readers can consult, for example, the book by Amari (1985).

4.3. *Pullback to* $\Gamma$. Using the embedding $\Lambda$ of $\Gamma$ in $\mathcal{M}(X)$, we can pull the metric on $\mathcal{M}(X)$ back to $\Gamma$. The definition of the pullback of the metric acting on two tangent vectors $u$ and $v$ in $T_\gamma \Gamma$ is as before

$$\mathbf{h}_\Lambda(u,v) = (\Lambda^*\mathbf{g})_\gamma(u,v) = \mathbf{g}_{\Lambda(\gamma)}(\Lambda_*(u)\Lambda_*(v)),$$

where $\Lambda_* : T_\gamma \Gamma \to T_{\Lambda(\gamma)}\mathcal{M}(X)$ is the tangent (derivative) map. This expression is coordinate-invariant. If we wish to know the matrix elements of $\mathbf{h}_\Lambda = \Lambda^*\mathbf{g}$ in the basis determined by a system of coordinates, $\frac{\partial}{\partial \theta^i}$, on $\Gamma$, we must evaluate $\mathbf{h}_\Lambda$ on these basis elements. The result is

$$\mathbf{h}_{\Lambda,\gamma}\left(\frac{\partial}{\partial \theta^i}(\gamma), \frac{\partial}{\partial \theta^j}(\gamma)\right) = \int_X \Lambda_\theta \frac{1}{\Lambda_\theta} \frac{\partial \Lambda_\theta}{\partial \theta^i} \frac{1}{\Lambda_\theta} \frac{\partial \Lambda_\theta}{\partial \theta^j},$$

where we denote by $\Lambda_\theta$ the value of the model function at the point $\gamma$ with coordinates $\theta$. We thus find the known result that the components of the induced metric form the Fisher information matrix.

As described in Section 3, the coordinate-invariant measure on $\Gamma$ is then given by

$$\mathbf{U}_\Lambda = \star_{\mathbf{h}_\Lambda}\mathbb{I} = |\mathbf{h}_\Lambda|_\theta^{1/2}\, d^m\theta.$$

4.4. *MAP estimates.* MAP estimation is now simply a question of using (3.1) with $\mathbf{Q}$ equal to the posterior measure from Bayes' theorem, and $\mathbf{U}$ equal to $\mathbf{U}_\Lambda$.

Note that the introduction of a prior probability prevents the estimate from being invariant under replacement of $\Lambda$ by $\Lambda\varepsilon$. The solution to this problem is the following. The prior probability is assigned to one member of the equivalence class $\{\Lambda\varepsilon\}$ based on knowledge of the parameters that is independent of current data. It can then be pushed forward to other copies of $\Gamma$ using $\varepsilon^{-1}$. Note that this violates condition 2 as it should, but that it does not violate condition 1.

In cases of "complete ignorance" of the value of $\gamma$, Jeffreys' prior is often used as the prior probability measure. In this case, the prior measure and the underlying measure cancel in the invariant MAP estimate, leaving only the model function. In cases of "complete ignorance," then MAP estimation reduces to maximum likelihood estimation *regardless* of the nature of Jeffreys' prior. (Note that the posterior probability measure still contains Jeffreys' prior; it is in the MAP estimate itself that it disappears.)

Thus while traditional MAP estimates of the variance of a Gaussian measure, for example, vary with the parameterization, the invariant MAP estimate will produce the maximum likelihood result in every case. The data space is $X = \mathbb{R}^n$, corresponding to $n$ independent experiments, and the model function is a Gaussian family of product measures on $n$, for the sake of argument with zero mean. The parameter space $\Gamma$ is isomorphic to $\mathbb{R}^+$: we use coordinates $\sigma \in \mathbb{R}$ on this space, where $\sigma$ is the standard deviation. The model function $\Lambda$ is then given by

$$\Lambda_\sigma = d^n x \, (2\pi\sigma^2)^{-n/2} \exp\left\{-\frac{(x,x)}{2\sigma^2}\right\},$$

where $(\cdot, \cdot)$ denotes the Euclidean inner product on $\mathbb{R}^n$. Derivation of the Fisher information then shows that the inner product between tangent vectors $u$ and $v$ in $T_\gamma\Gamma$, where the point $\gamma$ has coordinate $\sigma$, is

$$(4.1) \qquad \mathbf{h}_\Lambda(u,v) = \frac{2n}{\sigma^2} u^\sigma v^\sigma,$$

where the superscript $\sigma$ denotes the component with respect to the coordinate basis $\frac{\partial}{\partial\sigma}$. The induced measure is thus proportional to $d\sigma/\sigma$, the well-known Jeffreys' prior. Let us now consider the parameterization $v = \sigma^\alpha$, for $\alpha \in \mathbb{N}$. Jeffreys' prior is equal to $dv/v$ for all $\alpha \neq 0$. The traditional MAP estimates derived from these different parameterizations are

$$\hat{v}^{2/\alpha} = \frac{(x,x)}{n+\alpha},$$

where we have raised the estimate of $v$ to the power of $2/\alpha$ to make it equivalent to an estimate of $\sigma^2$. The problem of lack of invariance comes sharply into focus in this example. Which estimate of $\sigma$ is to be used?

On the other hand, the invariant MAP estimate is

$$\hat{v}^{2/\alpha} = \frac{(x,x)}{n}$$

for all $\alpha$.

4.5. *MMSD estimates.* In Section 3.2 we defined a coordinate-invariant version of the mean squared error estimate, which we called the MMSD estimate. Having defined a metric on $\Gamma$ above, we can now use it to calculate distances in $\Gamma$, and hence to define the MMSD estimate. In general, this is a difficult task that is not tractable analytically, although approximations may be available. In simple examples, however, one can compute the distance function $d(\gamma, \gamma')$ analytically. We give an example in Section 4.5.1.

4.5.1. *MMSD estimate of variance.* Consider the same example as above, of the estimation of the variance of a zero mean Gaussian measure.

From (4.1), the infinitesimal distance $ds$ between the points with coordinates $\sigma$ and $\sigma + d\sigma$ is given by

$$ds^2 = \frac{2n}{\sigma^2}\,d\sigma^2.$$

This is easily integrated to give the distance between two points with coordinates $\sigma_0$ and $\sigma_1$ (assume $\sigma_1 > \sigma_0$):

$$d(\sigma_0, \sigma_1) = \sqrt{2n}\ln\left(\frac{\sigma_1}{\sigma_0}\right).$$

The MMSD estimate of $\sigma$ is therefore given by considering the following mean loss under the posterior measure $\mathbf{Q}$ for $\sigma$:

$$L(\sigma) = \sqrt{2n}\int_0^\infty d\sigma' Q(\sigma')(\ln\sigma - \ln\sigma')^2.$$

Differentiation with respect to $\sigma$ then shows that the minimum squared distance estimate of $\sigma$, $\hat{\sigma}$, is given by

$$\hat{\sigma} = \exp E_{\mathbf{Q}}[\ln\sigma],$$

where $E_{\mathbf{Q}}[\cdot]$ indicates expectation using the measure $\mathbf{Q}$. Note that $E_{\mathbf{Q}}[\ln\sigma] \neq \ln E_{\mathbf{Q}}[\sigma]$ in general and that therefore the estimate is not simply the mean of $\sigma$ as would have been obtained by assuming a Euclidean metric.

The mean of $\ln\sigma$ can be calculated in the case that the prior on $\sigma$ is taken to be Jeffreys' prior. It is given in terms of coordinates by

$$E_{\mathbf{Q}}[\ln\sigma] = \tfrac{1}{2}[\ln(\tfrac{1}{2}(x,x)) - \psi(\tfrac{1}{2}n)],$$

where $\psi$ is the function

$$\psi(z) = \frac{d}{dz}\ln\Gamma(z)$$

and $\Gamma$ is the Gamma function $\Gamma(z) = \int_0^\infty dt\, t^{z-1}e^{-t}$. Thus

$$\hat{\sigma} = \sqrt{\frac{(x,x)}{2}}e^{-(1/2)\psi(n/2)}.$$

For large $z$, $\psi(z) \cong \ln(z)$, so that the estimate becomes

$$\hat{\sigma}_{\mathrm{cl}} = \sqrt{\frac{(x,x)}{2}}e^{-(1/2)\ln(n/2)} = \sqrt{\frac{(x,x)}{n}},$$

the classical result. To the next order, $\psi(z) \cong \ln(z) - \frac{1}{2z}$. This introduces corrections to the classical result:

$$\hat{\sigma} = e^{1/(2n)}\hat{\sigma}_{\mathrm{cl}}.$$

This formula is valid within about 10% down to $n = 1$, at which point the invariant result is bigger than the classical result by a factor of 1.9.

4.5.2. *General case in one dimension.* The form of the above estimate is quite general in the one-dimensional case. Consider that we have derived the metric on $\Gamma$, $\mathbf{h}$. The distance between two points $\gamma_0$ and $\gamma_1$ is then given according to the general discussion in Section 3. In a general coordinate system, $\theta$, this can be written

$$d(\gamma, \gamma') = \int_{t_0}^{t_1} dt \left( h(\pi(t)) \left( \frac{d\pi}{dt}(t) \right)^2 \right)^{1/2} = \int_{\theta_0}^{\theta_1} d\theta \, h^{1/2}(\theta),$$

where $\pi(t_{0,1}) = \gamma_{0,1}$, $\theta_{0,1} = \theta(\gamma_{0,1})$ and $h$ is the (single) component of the metric $\mathbf{h}$ in the $\theta$ coordinate system. Note that there is no need for a minimization in one dimension. All paths with the same endpoints belong to the same equivalence class under the action of (boundary- and orientation-preserving) diffeomorphisms of $I$. Now let $H$ be the inverse derivative of $h^{1/2}$. The (signed) distance between the two points is now $d(\theta_1, \theta_0) = H(\theta_1) - H(\theta_0)$. Including this in (3.3), differentiating $L$ and equating to zero then gives the result that

$$H(\hat{\theta}) = E_{\mathbf{Q}}[H],$$

and thus that

$$\hat{\theta} = H^{-1} E_{\mathbf{Q}}[H].$$

In more than one dimension, of course, the problem is a great deal more complicated, since there is an infinity of equivalence classes, and the minimization means solving a partial differential equation for the geodesics.

**5. Discussion and related work.** There is a significant amount of work on the geometry of probability measure spaces from the point of view of classical statistics; Murray and Rice (1993) and Kass and Vos (1997) provide recent treatments. As interesting as this work is, it has focused on asymptotics and other issues of importance to classical statistics, while the Bayesian approach using prior and posterior probabilities and loss functions has largely been ignored. As a consequence, it is not directly relevant to the problem posed in this paper. For example, Murray and Rice (1993) assert that the Riemannian distance is not of statistical significance, although they give no arguments, and that the mean in a manifold cannot be calculated; all that is possible is an analysis of the way in which the value of the mean, calculated in coordinates, changes with the coordinates. As we have seen, however, the Riemannian metric precisely allows the definition of a natural, coordinate-invariant generalization of the mean.

The pulled back metric defined in Section 4.3 was first introduced by Rao (1945), but it was the work of Amari (1985) that brought these ideas to prominence. Amari (1985) introduced, in addition to the metric, a family of connections on $\Gamma$, one of which was the metric connection compatible with the metric. The nonmetric connections, however, cannot be used to define the structures necessary for invariant Bayesian estimation as described here. Efron and Hinkley (1978) and Barndorff-Nielsen (1987) introduce "observed" geometric structures, but again these do not enable the definition of invariant estimates satisfying the two criteria in this paper. For example, the observed Fisher information metric of Efron and Hinkley (1978) is not a tensor, and thus violates the first criterion. In addition, it requires the definition of an underlying measure on the data space $X$; estimation is not invariant to this choice. Critchley, Marriott and Salmon (1994) develop "preferred point geometry" to try to ameliorate the lack of naturality they perceive in previous geometric approaches to statistics. The "preferred point metric" they define is, however, not invariant to diffeomorphisms, precisely because there is a preferred point. It thus violates the second criterion.

There is, from a Bayesian point of view, a more general objection to the asymmetric or preferred point structures (many of which also violate the triangle inequality) used in much of the above work. This objection is essentially the same as the original motivation for introducing them, which is the notion that there is a "true distribution" that must be treated differently, and related problems, for example, the worry that this distribution might not lie in the image of $\Gamma$. This notion does not exist, and indeed does not make sense, in a Bayesian approach. This can be seen by using, for example, a preferred point metric in the formula for the posterior density, (3.1). The preferred point is undefined, yet if it is taken to be the argument to the posterior density, seemingly the only reasonable choice, then the "preferred point" vanishes and we are back to the Riemannian metric described herein. Thus the raison d'être of these more complex structures disappears.

From another direction, Pennec (1999) develops some basic statistical tools for Riemannian manifolds, and applies these ideas in various ways to problems in computer vision. The approach is not Bayesian, however, and in particular the choice of a metric and the relation with estimation problems, including the use of the metric measure as an underlying measure for MAP estimation, are not considered.

MML inference was developed by Wallace and Boulton (1968) and Wallace and Freeman (1987). A discussion of its relationship with the standard Bayesian approach and of its invariance properties can be found in the above papers and in the paper by Oliver and Baxter (1995). The literature on MML inference frequently cites the invariance of MML estimates as one reason to prefer them to MAP estimates. The above analysis shows that this is not a special property of MML estimates, or a deep problem with MAP estimates. Indeed, the issue is not one of MAP estimation per se. Lack of invariance is a consequence of not describing the quantities of interest in $\Gamma$ in a coordinate-invariant, and hence meaningful, way. To do this, one must recognize that a metric is lurking in the definition of both MAP and MMSE estimates, and indeed in any useful discussion of $\Gamma$, and that making it explicit is a necessary condition for meaningful definitions in arbitrary coordinate systems. Once done, the definition of coordinate-invariant estimates is an immediate consequence of the geometry. Although (3.1) with the pulled back metric as underlying measure is formally the same as that for MML estimates, unlike MML methods, no discretization of $\Gamma$ is needed, and no approximations are made. In fact, the above derivation throws light on the procedure used in deriving MML estimates, which from this point of view appears to be a roundabout way of defining an underlying measure by first discretizing the manifold and then considering the volume of each cell.

The fact that we are discussing the geometry of $\Gamma$ and not a particular form of estimate means that the analysis presented here is more general than MML, however. By recognizing the necessity of an explicit metric on $\Gamma$ for inference, the way is open for the definition of coordinate-invariant loss functions of many different types. Here we have given the example of a coordinate-invariant MMSE estimate, the MMSD estimate, but whenever defining a loss function on a parameter space, the issues described here must be taken into account.

5.1. *Discussion of choice of metric.* In Section 4 we came to the conclusion that the only choice of metric that satisfies the two conditions mentioned at the beginning of that section is the metric induced by pullback from $\mathcal{M}(X)$. To recap: the metric and its associated underlying measure should not introduce information about $\Gamma$. Such information should be contained in one of two sources: the correspondence between points of $\Gamma$ and points of $\mathcal{M}(X)$, and the prior measure. The first leads to the idea that the metric on diffeomorphically related copies of $\Gamma$ should be related by pullback, while the second eliminates the possibility of choosing a metric on one fixed copy of $\Gamma$ and then pulling it back to the other copies, since this implies that we must be able to assign a value of the metric to particular points in $\Gamma$ a priori, which in turn implies that we must know something about the identity of these points beyond the information contained in the prior. Hence the result given.

Note that this argument is somewhat different from that normally used for Jeffreys' prior, or rather is a clarification and a refinement of that argument, which essentially boils down to proving that this prior is invariant under "reparameterizations." First, the emphasis is on the metric as providing $\Gamma$ with geometry, and not on the measure, which is a derived quantity. Second, coordinate invariance is not an issue: the abstract way in which the geometry is described does not rely on a particular choice of coordinate system. Equation (3.1), for example, is coordinate-invariant for any choice of metric. Instead the emphasis is on diffeomorphism invariance: our results should not depend on which copy of $\Gamma$ we use, since this merely "shuffles" the points of $\Gamma$ without changing their correspondence with points of $\mathcal{M}(X)$.

The use of the underlying measure of the pulled back metric does not commit us to using Jeffreys' prior as a noninformative prior. Thus the large amount of previous work [Bernardo (1979) and Kass and Wasserman (1996)] on the choice of such priors, fascinating though it is, is not directly relevant to our discussion here. Note in particular that the problems associated with Jeffreys' prior do not appear when we are talking about an underlying measure. Normalization is not necessary since the underlying measure is not a probability measure. Second, the procedure advocated here suggests that we should first eliminate nuisance parameters using whatever prior information we possess, to obtain a likelihood on the parameter of interest, and only then derive the metric by pullback. Thus the various "paradoxes" associated with the noncommutativity of the derivation of Jeffreys' prior and marginalization do not arise.

Our argument for the metric and underlying measure on $\Gamma$ does not depend on group-theoretic considerations. Nevertheless, the metric is compatible with these considerations, as is Jeffreys' prior, because of the following simple argument. Let $X$ be a manifold with metric $\mathbf{h}$, and let $Y$ be embedded in $X$ by $f$. Suppose we have two group actions $\beta_X : G \times f(Y) \to f(Y)$ and $\beta_Y : G \times Y \to Y$. Note that the group action on $X$ need only be defined for the image of $Y$; it may, for example, be induced by the group action on $Y$ itself. If we have

$$
\begin{array}{ccc}
Y & \xrightarrow{\ f\ } & f(Y) \\
{\scriptstyle \beta_Y(g)} \big\uparrow & & \big\uparrow {\scriptstyle \beta_X(g)} \\
Y & \xrightarrow[\ f\ ]{} & f(Y)
\end{array}
$$

then, if $G$ acts by isometries on $X$, endowing $Y$ with the metric $f^*\mathbf{h}$ ensures that $f$ is an isometry also. Therefore, $G$ must act by isometries on $Y$. If $Y$ is $G$ itself, this ensures that the underlying measure induced by the metric $f^*\mathbf{h}$ is a Haar measure.

Finally, an information-theoretic intuition is interesting. In computing the MAP estimate, it is equivalent to maximize the logarithm of (3.2). Naturally the logarithm consists of the difference of two terms: the logarithm of the posterior density and the logarithm of the underlying density. The role of the underlying density is the following. The information that we possess should presumably be that amount of information that we possess beyond "ignorance." If our expression for "ignorance" does not possess the value "zero" (i.e., the identity) in the algebra in which we add and subtract information, then the information that we possess beyond "ignorance" should be the difference between the algebraic element representing our knowledge, and the algebraic element representing "ignorance." In view of the "noninformative" nature of the underlying measure that we are using, the MAP estimate can thus consistently be thought of as finding that point in $\Gamma$ with maximum information.

This intuition, and the invariant nature of the underlying measure, suggest that this measure should be the reference measure for the maximum entropy approach to generating prior measures on manifolds. This is a subject for further research.

## APPENDIX: FORMS

We provide a short introduction to the language of forms. A good reference is the book by Choquet-Bruhat, DeWitt-Morette and Dillard-Bleick (1977). Briefly, differential forms are antisymmetric, multilinear functionals on products of vector spaces. For manifolds they are defined pointwise on the tangent space at each point and then required to satisfy smoothness properties. They also allow a beautiful theory of integration on manifolds, and in this capacity they are thought of as *co-chains*, linear functionals on the vector space of chains in a manifold. Their advantages are great concision and uniformity of notation; independence of basis or coordinates; manifest invariance to diffeomorphisms and other transformations; and generality. In bringing together integration and geometry in one notation, they are ideal for our discussion.

We are given a manifold $\Gamma$. From here, we can define the tangent space at each point, $T_\gamma \Gamma$, using a number of approaches. The result is intuitively clear, however, so we will not go into detail. We can bring all the tangent spaces together in the *tangent bundle*, $T\Gamma$. This is another manifold, each point of which can be thought of as a pair: a point $\gamma$ in $\Gamma$ and a vector in $T_\gamma \Gamma$. There is a canonical projection from $T\Gamma$ to $\Gamma$ supplied by forgetting the tangent vector. At each point $\gamma$, the tangent space $T_\gamma \Gamma$ has a dual space, $T_\gamma^* \Gamma$, the space of linear maps from $T_\gamma \Gamma$ to $\mathbb{R}$. These can be combined to form the co-tangent bundle, $T^* \Gamma$. A *vector field* is a *section* of the tangent bundle: a map from $\Gamma$ to $T\Gamma$ whose left inverse is the canonical projection.

We can also form product bundles, in which the "extra space" at each point $\gamma$ is the product of copies of the tangent space; thus each point in $T^p \Gamma$ can be thought of as a point $\gamma$ and an element of $\bigotimes^p T_\gamma \Gamma$. Now at each point we can define higher dual spaces: $T_\gamma^{*p} \Gamma = \bigotimes^p T_\gamma^* \Gamma$ is the space of multilinear functions on $\times^p T_\gamma \Gamma$. In particular, we can restrict attention to the antisymmetric linear functions: those that change sign under the interchange of any two arguments. These are antisymmetric tensor products of the co-tangent space, denoted $\bigwedge^p T_\gamma^* \Gamma$. Their combination into a bundle is denoted $\bigwedge^p T^* \Gamma$. A section of $\bigwedge^p T^* \Gamma$ defines, for each point $\gamma$, an element of $\bigwedge^p T_\gamma^* \Gamma$. Sections of $\bigwedge^p T^* \Gamma$ are known as *forms*, and $p$ is the *degree* of the form. We denote the space of $p$-forms $\Lambda^p \Gamma$. Forms of degree $p$ and $q$ can be multiplied to give forms of degree $p + q$. Because the product of co-tangent spaces is antisymmetric, all forms of degree higher than $m$, the dimensionality of the manifold, are zero. 0-forms are functions on $\Gamma$.

In order to express vectors and forms more easily, it is convenient to introduce bases for the various spaces. This is easily done using a coordinate system $\theta : \Gamma \to \mathbb{R}^m$. A basis for $T_\gamma \Gamma$ is then the set of $\frac{\partial}{\partial \theta^j}(\gamma)$. The dual basis for $T_\gamma^* \Gamma$ is then the set of $d\theta^i(\gamma)$, which acts on the basis of $T_\gamma \Gamma$ as

$$d\theta^i(\gamma) \left( \frac{\partial}{\partial \theta^j}(\gamma) \right) = \delta^i_j.$$

Taking the collection of these bases all over $\Gamma$, we have bases for the spaces of vector fields and 1-forms. Now we can form bases for the various power bundles. For example, a basis for the space of 2-forms is given by the set $d\theta^i(\gamma) \wedge d\theta^j(\gamma)$, where $\wedge$ denotes the antisymmetric product. We will denote the basis element $d\theta^1(\gamma) \wedge \cdots \wedge d\theta^m(\gamma)$ of the space of $m$-forms (there is only one—if the indices are not different, antisymmetry of the product means the result is zero) by $d^m \theta(\gamma)$. The sign of this basis element (or in other words, the order of the factors of $d\theta^i$ that it contains) defines an *orientation* on the manifold, in the sense that a basis for the tangent spaces, when acted upon by the form, will give either a positive or negative result depending on its orientation in the traditional sense of right- and left-handed coordinate systems. Given an orientation in this sense, a basis for the tangent spaces is either *oriented* or not. Not all manifolds admit a global orientation. We consider only orientable manifolds.

Given another manifold $Y$, and a map $\Lambda : Y \to \Gamma$, we define the *tangent map* or *derivative map* at a point $y \in Y$, $\Lambda_* : T_y Y \to T_{\Lambda(y)} \Gamma$ as follows. A point $(y, u) \in TY$ is taken to $(\Lambda(y), \Lambda_* u) \in T\Gamma$, where, in terms of coordinates $\theta^i$ on $\Gamma$ and $\phi^\alpha$ on $Y$, in which $u = u^\alpha \frac{\partial}{\partial \phi^\alpha}$, we have

$$\Lambda_* u = (\Lambda_* u)^i \frac{\partial}{\partial \theta^i} = u^\alpha \frac{\partial \Lambda^i}{\partial \phi^\alpha} \frac{\partial}{\partial \theta^i},$$

where $\Lambda^i = \theta^i(\Lambda)$. We also introduce the convention that repeated indices, one up, one down, are summed over.

Using this map, we can define the *pullback* $\Lambda^* \mathbf{A}$ of a form $\mathbf{A} \in \Lambda^p \Gamma$ (or in fact of any member of a power of a co-tangent space, whether antisymmetric or not) as

$$\Lambda^* \mathbf{A}_y(u, v, \dots) = \mathbf{A}_{\Lambda(y)}(\Lambda_* u, \Lambda_* v, \dots).$$

Thus the action of a pulled back form on tangent vectors is defined by the action of the original form on the tangent vectors pushed forward by the tangent map.

As well as antisymmetric products of co-tangent spaces, we can form symmetric products. If at each point $\gamma$ we form the space of symmetric, bilinear functions on $T_\gamma\Gamma \times T_\gamma\Gamma$, which we will denote $T_\gamma^*\Gamma \vee T_\gamma^*\Gamma$, we can again form a product bundle $T^*\Gamma \vee T^*\Gamma$. A *metric* $\mathbf{h}$ on $\Gamma$ is a positive (semi-)definite section of this bundle: to each point $\gamma$ it assigns a positive (semi-)definite element of $T_\gamma^*\Gamma \vee T_\gamma^*\Gamma$, or in other words, an inner product on $T_\gamma\Gamma$.

In a particular coordinate basis $\frac{\partial}{\partial\theta^i}(\gamma)$, the metric has components given by

$$h_{\gamma,ij} = \mathbf{h}_\gamma\left(\frac{\partial}{\partial\theta^i}(\gamma), \frac{\partial}{\partial\theta^j}(\gamma)\right).$$

The matrix elements of the metric at each point $\gamma$ possess a determinant, which we will write $|\mathbf{h}|_\theta(\theta(\gamma))$.

Using the metric $\mathbf{h}$, we can define a canonical isomorphism, the *Hodge star* $\star_{\mathbf{h}}$, between $\Lambda^p\Gamma$ and $\Lambda^{m-p}\Gamma$. We show here its action for $p = 0$ and $p = m$ only, since that is all we will need. We choose coordinates $\theta^i$ (nothing will depend on this choice). Let $f$ be a 0-form, and let $\mathbf{A} = A\,d^m\theta$ be an $m$-form ($A$ is a function—the component of $\mathbf{A}$ in the basis $d^m\theta$). Then we have

(A.1)
$$\begin{aligned}\star_{\mathbf{h}} f &= |\mathbf{h}|^{1/2} f \, d^m\theta, \\ \star_{\mathbf{h}} \mathbf{A} &= |\mathbf{h}|^{-1/2} A,\end{aligned}$$

where we have suppressed arguments and reference to the coordinate system in the definition of the determinant for clarity.

The Hodge star can be used to define an inner product on each $\Lambda^p\Gamma$. Since $\star_{\mathbf{h}}\mathbf{A}$ is an $(m-p)$-form if $\mathbf{A}$ is a $p$-form, the quantity $\mathbf{A}\star_{\mathbf{h}}\mathbf{B}$ for two $p$-forms is an $m$-form, and can be integrated on $\Gamma$:

(A.2)
$$\langle\!\langle\mathbf{A}, \mathbf{B}\rangle\!\rangle = \int_\Gamma \mathbf{A}\star_{\mathbf{h}}\mathbf{B}.$$

We can define *positive $m$-forms* as those whose action on oriented bases produces a positive result. It is equivalent to say that their dual under the action of the Hodge star is a positive function. A *probability $m$-form* is a positive $m$-form whose integral over $\Gamma$ is equal to 1. We can divide $m$-forms by positive $m$-forms. For an $m$-form $\mathbf{A}$ and a positive $m$-form $\mathbf{B}$, the value of $\frac{\mathbf{A}}{\mathbf{B}}$ is that unique function $f$ such that $\mathbf{A} = f\mathbf{B}$. This division is the analogue of the Radon–Nikodym derivative for forms.

On an $m$-dimensional manifold, $m$-forms can be integrated in the way that the notation suggests. For an $m$-form $\mathbf{A} = A\,d^m\theta$, we have that

$$\int_{\Omega\subset\Gamma}\mathbf{A} = \int_{\theta(\Omega)} A(\theta)\,d^m\theta,$$

where we have used the same symbol $A$ for the function and its expression in terms of coordinates.

To integrate a $p$-form $\mathbf{A}$ over a $p$-dimensional submanifold embedded in $\Gamma$, $Y \overset{\Lambda}{\rightarrowtail} \Gamma$, one first pulls the form back to the embedded manifold and then integrates:

$$\int_{\Lambda(Y)} \mathbf{A} = \int_Y \Lambda^* \mathbf{A}.$$

These definitions highlight the second way of interpreting forms: as *co-chains*. A $p$-chain in $\Gamma$ is (roughly speaking) a linear combination of $p$-dimensional rectangles embedded in the manifold. The space of linear functions on the space of $p$-chains (the co-chains) can be identified with $\Lambda^p \Gamma$.

We will have cause to integrate a function $f$ over a $p$-dimensional submanifold $Y \overset{\Lambda}{\rightarrowtail} \Gamma$ of $\Gamma$. This is slightly different from the case of integrating a $p$-form. One first pulls the function back to $Y$ and then uses a metric on $Y$ to convert the function into a $p$-form that can be integrated over $Y$:

$$\int_{\Lambda(Y)} f = \int_Y \star_{\mathbf{h}} \Lambda^* f,$$

where by definition $(\Lambda^* f)(y) = f(\Lambda(y))$, and $\mathbf{h}$ is a metric on $Y$.

However, since we are interested in the submanifold in $\Gamma$ and not $Y$ itself, we are really considering an equivalence class of embeddings $\{f\varepsilon\}$, where $\varepsilon : Y \to Y$ is a diffeomorphism, with the same image. The result of our integration should be independent of the representative in this equivalence class, and this means that the metric on $Y$ must vary with the representative. If no representative is distinguished, the only way to achieve this invariance is to pull back a metric $\mathbf{g}$ on $\Gamma$ to $Y$, and use this metric to define the Hodge star:

$$\int_{\Lambda(Y)} f = \int_Y \star_{\Lambda^* \mathbf{g}} \Lambda^* f.$$

**Acknowledgment.**  The author would like to thank both referees for their valuable comments on the paper.

ARIANA (JOINT RESEARCH PROJECT INRIA/I3S)
INRIA
2004 ROUTE DES LUCIOLES
B.P. 93
06902 SOPHIA ANTIPOLIS CEDEX
FRANCE
E-MAIL: Ian.Jermyn@sophia.inria.fr